\documentclass{article}

\usepackage{amsthm,amsmath,amssymb}

\newtheorem{theorem}{Theorem}[section]

\newtheorem{proposition}[theorem]{Proposition}
\newtheorem{corollary}[theorem]{Corollary}
\newtheorem{lemma}[theorem]{Lemma}

\def\R{{\mathbb R}}
\def\N{{\mathbb N}}

\newcommand{\remin}{\mathop{-\!\!\!\!\!\hspace*{1mm}\raisebox{0.5mm}{$
\cdot$}}\nolimits}

\begin{document}

\title{The approximate fixed point property in product spaces}
\author{U. Kohlenbach$^{1}$, L. Leu\c stean$^{1,2}$\\[0.2cm]
\footnotesize ${}^1$ Department of Mathematics, Darmstadt University of Technology,\\
\footnotesize Schlossgartenstrasse 7, 64289 Darmstadt, Germany\\[0.1cm]
\footnotesize${}^2$ Institute of Mathematics "Simion Stoilow'' of the Romanian Academy, \\
\footnotesize Calea Grivi\c tei 21, P.O. Box 1-462, Bucharest, Romania\\[0.1cm]
\footnotesize E-mails: kohlenbach,leustean@mathematik.tu-darmstadt.de
}

\date{}
\maketitle

\begin{abstract}
In this paper we generalize to {\em unbounded} convex subsets $C$ of {\em hyperbolic} spaces results obtained by W.A. Kirk and  R. Esp\'{\i}nola on approximate fixed points of nonexpansive mappings in product spaces $(C\times M)_\infty$, where $M$ is a metric space and $C$ is a nonempty, convex, closed and bounded subset of a normed or a CAT(0)-space. We extend the results further, to families $(C_u)_{u\in M}$ of unbounded convex subsets of a hyperbolic space. The key ingredient in obtaining these generalizations  is a uniform quantitative version of a theorem due to Borwein, Reich and Shafrir, obtained by the authors  in a previous paper using techniques from mathematical logic. Inspired by that, we introduce in the last section  the notion of {\em uniform approximate fixed point property} for sets $C$ and classes of self-mappings of $C$.  The paper ends with an open problem.
\end{abstract}
\begin{tabular}{ll}
\footnotesize\noindent {\it Keywords}: &\footnotesize Proof mining, metric fixed point theory, nonexpansive functions, \\
&\footnotesize  approximate fixed points,  product spaces\\[0.1cm]
\noindent {\it MSC:\ } & \footnotesize 47H10,  47H09, 03F10
\end{tabular}

\section{Introduction}

This paper presents applications of a case study in the project of 
{\em proof mining}, 
by which we mean the logical analysis of mathematical proofs with the aim of 
extracting 
new numerically relevant information hidden in the proofs.

More specifically, we are concerned with the 
approximate fixed points for nonexpansive mappings in product spaces. 

The main tool used in this paper is a quantitative version (Theorem \ref{quant-BRS}) of a theorem due to 
\cite{Borwein/Reich/Shafrir} (see Theorem \ref{brs-th}) which was obtained in \cite{Kohlenbach/Leustean(03)} by logical analysis of the 
original proof of Theorem \ref{brs-th}.

Let $(X,\rho)$ be a metric space, and $C\subseteq X$ a nonempty subset. A mapping $T:C\to C$ is called {\em nonexpansive} if for all $x,y\in C$,
\[\rho(T(x), T(y))\le \rho(x,y).\]
The metric space $(X,\rho)$ is said to have the {\em approximate fixed point property (AFPP)} for nonexpansive mappings if each nonexpansive mapping $T:X\to X$
has an approximate fixed point sequence; that is, a sequence $(u_n)_{n\in\N}$ in $X$ for which $\displaystyle\lim_{n\to\infty} \rho(u_n, T(u_n))=0$. It is easy to see that this is equivalent with  
\[r_X(T):=\inf\{\rho(x,T(x))\mid x\in X\}=0.\]
If $(X,\rho)$ and $(Y,d)$ are metric spaces, then the metric $d_\infty$ on $X\times Y$ is defined in the usual way:
\[ d_\infty ((x,u),(y,v))=\max\{\rho(x,y), d(u,v)\}\]
for $(x,u), (y,v)\in X\times Y$. We denote by $(X\times Y)_\infty$ the metric space thus obtained.

The following theorem was proved first by Esp\' inola and Kirk \cite{Espinola/Kirk} for normed spaces and then proved by Kirk 
\cite[Theorem 25, Remark 26]{Kirk(04)} for all CAT(0)-spaces. 
\begin{theorem}\label{th-kirk-CAT}
Assume that $X$ is a normed space or a CAT(0)-space,  
$C\subseteq X$ is  a nonempty, convex, closed and {\em bounded} subset of $X$. 
If $(M,d)$ is a metric space with the AFPP for nonexpansive mappings, then
\[H:=(C\times M)_\infty\]
has the AFPP for nonexpansive mappings.
\end{theorem}
Kirk's proof in \cite{Kirk(04)} actually also holds for the class of {\em 
hyperbolic spaces} as introduced by the first author in \cite{Kohlenbach(05a)}: 

A {\em  hyperbolic space}\footnote
{ See 
\cite{Kohlenbach(05a)} for discussion of this and related notions.} 
is a triple $(X,\rho,W)$ where
$(X,\rho)$ is a metric space and $W:X\times X\times [0,1]\to X$ is such that 
\begin{eqnarray*}
(W1) & \rho(z,W(x,y,\lambda))\le (1-\lambda)\rho(z,x)+\lambda \rho(z,y),\\
(W2) & \rho(W(x,y,\lambda),W(x,y,\tilde{\lambda}))=|\lambda-\tilde{\lambda}|\cdot 
\rho(x,y),\\
(W3) & W(x,y,\lambda)=W(y,x,1-\lambda),\\
(W4) & \,\,\,\rho(W(x,z,\lambda),W(y,w,\lambda)) \le (1-\lambda)\rho(x,y)+\lambda
\rho(z,w).
\end {eqnarray*}
If only axiom $(W1)$ is assumed this structure is a convex metric space in the sense of Takahashi \cite{Takahashi}. If $(W1)$-$(W3)$ are assumed, the notion is equivalent to Kirk's spaces of hyperbolic type \cite{Goebel/Kirk}. 
Axiom $(W4)$ is used e.g. in \cite{Kirk(82),Reich/Shafrir}. However, the concept of hyperbolic space as defined 
in \cite{Reich/Shafrir} is somewhat more restrictive than ours, as it is assumed that the space contains 
metric lines. This has the consequence that only CAT(0)-spaces which have the unique geodesic line extension 
property are included whereas in our definition all CAT(0)-spaces are covered. The notion of space of 
hyperbolic type would be too general though for our purposes as the proof of Theorem \ref{brs-th} uses 
(W4) in an essential way. 

The class of hyperbolic spaces contains all normed linear spaces and convex subsets thereof, but also the open unit ball in complex Hilbert spaces with the hyperbolic metric as well as Hadamard manifolds and CAT(0)-spaces in the sense of Gromov.

{\bf Remark:} Kirk actually claims that his proof would go through for the 
still larger class of spaces of hyperbolic type as defined in 
\cite{Goebel/Kirk}. However, the proof of Theorem 25 in \cite{Kirk(04)} 
uses axiom $(W4)$ for the proof of Lemma 3 (used for 
$\lambda=\frac{1}{2}$) which does not seem to hold in this setting in general. 
The proof does work, though, for the more restricted notion defined in 
\cite{Kirk(82)} (and called there space of hyperbolic type as well) which 
is the same notion as the one used in \cite{Reich/Shafrir}. But this includes 
only CAT(0)-spaces with the geodesic line extension property. Our more liberal 
notion avoids this and Kirk's proof, nevertheless, goes through.

If $x,y\in X$, and $\lambda\in[0,1]$ then we use the notation $(1-\lambda)x\oplus \lambda y$ for $W(x,y,\lambda)$. 

If $C\subseteq X$ is a nonempty convex subset of a hyperbolic space $(X,\rho,W)$, and $T:C\to C$ is nonexpansive, then for each sequence $(\lambda_n)_{n\in\N}$ in $[0,1]$, we can define the {\em Krasnoselski-Mann iteration} starting from $x\in C$ by
\begin{equation}
x_0:=x, \quad x_{n+1}:=(1-\lambda_n)x_n\oplus \lambda_n T(x_n). \label{KM-lambda-n-def}
\end{equation}

An important result in the fixed point theory for nonexpansive mappings is 
the following theorem, due to Borwein, Reich, and 
Shafrir\footnote{ In \cite{Borwein/Reich/Shafrir}, Borwein, Reich and Shafrir actually only state the result for their more restricted notion of hyperbolic space as defined in \cite{Kirk(82)}. 
However, from \cite{Kohlenbach/Leustean(03)} it follows that the result 
holds for our concept of hyperbolic space as well.} (generalizing earlier results of Ishikawa \cite{Ishikawa(76)},  Goebel and Kirk \cite{Goebel/Kirk}).

\begin{theorem}\cite{Borwein/Reich/Shafrir}\label{brs-th}
If $(\lambda_n)_{n\in\N}$ is a sequence in $[0,1]$ which is divergent in sum 
and bounded away from $1$ then for all $x\in C$, 
\[\displaystyle\lim_{n\to\infty} \rho(x_n,T(x_n))=r_C(T). \] 
\end{theorem}
In  \cite{Kohlenbach/Leustean(03)}, we obtained (even for the more general class of directionally nonexpansive functions as introduced in \cite{Kirk(00)}) the following quantitative version of Theorem \ref{brs-th} (which subsequently turned out to be an instance of a general logical metatheorem, see \cite{Gerhardy/Kohlenbach(05)}):

\begin{theorem}\label{quant-BRS}
Let $(X,\rho, W)$ be a hyperbolic space, $C\subseteq X$ a nonempty convex subset and $T:C\rightarrow C$ a nonexpansive mapping. Let $(\lambda_n)_{n\in\N}$ be a sequence in $[0,1)$, and $K\in\N, \alpha:\N\to\N$ be such that for all $n\in\N$,
\begin{equation}
\lambda_n \le 1-\frac{1}{K}, \text{ and } n\le\sum\limits^{\alpha(n)}_{i=0}\lambda_i.\label{hyp-K-alpha}
\end{equation}
Let $x,x^*\in C$ and $b>0$ be such that
\[ \rho(x, T(x))\le b \wedge \rho(x,x^*)\le b.\]
Then the  following holds
\[ \forall \varepsilon >0 \forall n\ge
h(\varepsilon,b,K,\alpha)\big(\rho(x_n,T(x_n))\leq \rho(x^*,T(x^*))
+\varepsilon\big), \] 
where\footnote{$n\, \remin 1=\max\{0,n-1\}.$}
\[
\begin{array}{l}
h(\varepsilon,b,K,\alpha)=\widehat{\alpha}(\lceil
2b\cdot \exp(K(M+1)) \rceil\remin\, 1,M), \  \mbox{with} \\
M\in \N \ \mbox{is such that} \ \displaystyle M\ge \frac{1+2b}{\varepsilon},\\
\widehat{\alpha}(0,n):=\tilde{\alpha}(0,n), \,\,
\widehat{\alpha}(i+1,n):=
\tilde{\alpha}(\widehat{\alpha}(i,n),n), \\
\tilde{\alpha}(i,n):=i+\alpha^+(i,n),\,\, \alpha^+(i,n):=\max\{\alpha'(j,n):j\le i\},\\
 \alpha'(i,n):=\alpha(n+i)-i+1, \quad (i,n\in\N).
\end{array}\]
\end{theorem}
Using Theorem \ref{quant-BRS}, the following quantitative version of a theorem 
due to Ishikawa \cite{Ishikawa(76)} (for normed spaces) and 
Goebel and Kirk \cite{Goebel/Kirk} (for  spaces of hyperbolic type) was 
obtained in \cite{Kohlenbach/Leustean(03)} (again the general form 
of this quantitative version is guaranteed by a logical metatheorem, 
see \cite{Gerhardy/Kohlenbach(05)}): 

\begin{theorem}\label{quant-Ishikawa}\cite{Kohlenbach/Leustean(03)}
Let $(X,\rho, W), C, (\lambda_n), K, \alpha$ be as in the previous theorem. 
Let $b>0, x,x^*\in C$ be such that 
\begin{equation}
\rho(x,x^*)\le b\wedge\forall m,n\in\N(\rho(x_n^*, x_m^*)\le b).
\end{equation}
 Then the  following holds
\[ \forall \varepsilon >0 \forall n\ge
\tilde{h}(\varepsilon,b,K,\alpha)\ \big(\rho(x_n,T(x_n))\leq\varepsilon\big), \] 
where 
\[
\begin{array}{l}
\tilde{h}(\varepsilon,b,K,\alpha)=\widehat{\alpha}(\lceil
12b\cdot \exp(K(M+1)) \rceil\remin\, 1,M), \  \mbox{with} \\
M\in \N \ \mbox{is such that} \ \displaystyle M\ge \frac{1+6b}{\varepsilon}, \\
\text{and } \widehat{\alpha} \text{~as before}.
\end{array}
\]
\end{theorem}
The main significance of the bounds in the previous theorems is that they depend on $x,x^*,T,C,X$ only via $b$ (and on $(\lambda_n)$ only via 
$K,\alpha$). In particular, if in 
Theorem \ref{quant-Ishikawa}, $C$ is assumed to have a bounded diameter, then 
the convergence $\rho(x_n,T(x_n))\to 0$ is uniform in $x$ and $T.$ This result 
was first obtained (ineffectively) in \cite{Goebel/Kirk} and used in 
\cite{Kirk(04)} to prove Theorem \ref{th-kirk-CAT} discussed above. 

In this paper, we shall use our strong uniformity results to approach (in the more general case of hyperbolic spaces) the following problem : 
what results can we obtain if we drop the hypothesis of $C$ being bounded from Theorem \ref{th-kirk-CAT}?

\section{Main results}

In the following, $(X,\rho,W)$ always is a hyperbolic space.

\subsection{The case of one convex subset $C$}

Let $C\subseteq X$ be a nonempty convex subset. We assume that $(M,d)$  is a metric space which has the AFPP for nonexpansive mappings.
Let $H:=(C\times M)_\infty$ and  $P_1:H\to C, P_2:H\to M$ be the coordinate projections. Finally, let $(\lambda_n)$ be a sequence in $[0,1].$

For each nonexpansive function $T:H\to H$ and for  each $u\in M$, let us define 
\[T_u:C\to C, \quad T_u(x)=(P_1\circ T)(x,u).\]
It is easy to see that $T_u$ is nonexpansive, so we can associate with $T_u$ the Krasnoselski-Mann iteration $(x_n^u)$  starting with an arbitrary $x\in C$.

In the sequel, $\delta:M\to C$ is a nonexpansive function, which {\em selects} for each $u\in M$ an element $\delta(u)\in C$. Trivial examples of such nonexpansive selection functions are the constant ones. 

For simplicity, we shall denote the Krasnoselski-Mann iteration starting from $\delta(u)$ and associated with $T_u$ by $(\delta(u))_n$.

For each $n\ge 0$, we define 
\begin{eqnarray*}
\varphi_n:M\to M, & \varphi_n(u)=(P_2\circ T)((\delta(u))_n,u).
\end{eqnarray*}

The following lemma collects some properties we will use in our proofs.

\begin{lemma}\label{useful-lemma}$\hfill$
\begin{enumerate}
\item[(i)] for all $u\in M, x,y\in C$, the sequences $(\rho(x_n^u, T_u(x_n^u)))$ and $(\rho(x_n^u, y_n^u))$  are nonincreasing;
\item[(ii)] for all $n\ge 0,\, u,v\in M$, and $x,y\in C$
\begin{equation}
\rho(x_n^u, y_n^v)\le \max\{\rho(x,y),d(u,v)\}=d_\infty((x,u),(y,v)); \label{ineq-x-u-y-v}
\end{equation}
\item[(iii)] for all $n\ge 0$, $\varphi_n$ is nonexpansive;
\item[(iv)] there is a sequence $(z_n)$ in $M$ such that for all $n\ge 0$,
\[\begin{array}{l}
\displaystyle d(z_n, \varphi_n(z_n))\le \frac{1}{n}, \text{and}\\
\displaystyle d_\infty(((\delta(z_n))_n, z_n), T((\delta(z_n))_n, z_n))\le \max\left\{\rho\big((\delta(z_n))_n, T_{z_n}((\delta(z_n))_n)\big), \frac{1}{n}\right\}.
\end{array}
\]
\end{enumerate}
\end{lemma}
\begin{proof}
\begin{enumerate}
\item[(i)] Apply, for example, \cite[Proposition 3.4, Lemma 3.8]{Kohlenbach/Leustean(03)}, whose proofs immediately generalize to our notion of hyperbolic space. 
\item[(ii)] The proof is by induction on $n$. The case $n=0$ is immediate. 
Assume that (\ref{ineq-x-u-y-v}) is true for $n$. We get
\begin{eqnarray*}
\rho(x_{n+1}^u, y_{n+1}^v) & = & \rho((1-\lambda_n)x_n^u
\oplus\lambda_nT_u(x_n^u),(1-\lambda_n)y_n^v\oplus\lambda_nT_v(y_n^v))\\
&\le  & (1-\lambda_n)\rho(x_n^u,y_n^v)+\lambda_n\rho(T_u(x_n^u),T_v(y_n^v)) 
\quad\mbox{\rm by (W4)} \\
&\le & (1-\lambda_n)\max\{\rho(x,y),d(u,v)\}+\\
&&+\lambda_n d_\infty(T(x_n^u,u),T(y_n^v,v))\\
&\le & (1-\lambda_n)\max\{\rho(x,y),d(u,v)\}+\\
&&+\lambda_n d_\infty((x_n^u,u),(y_n^v,v))\\
&= & (1-\lambda_n)\max\{\rho(x,y),d(u,v)\}+\\
&&+\lambda_n \max\{\rho(x_n^u,y_n^v),\,d(u,v)\}\\
& \le & (1-\lambda_n)\max\{\rho(x,y),d(u,v)\}+\\
&&+\lambda_n \max\{\rho(x,y),d(u,v)\}\\
&= & \max\{\rho(x,y),d(u,v)\}.
\end{eqnarray*}
\item[(iii)] We get that
\begin{eqnarray*}
d(\varphi_n(u),\varphi_n(v))&\le & d_\infty(T((\delta(u))_n,u),T((\delta(v))_n,v))\\
&\le &  d_\infty(((\delta(u))_n,u),((\delta(v))_n,v))\\
&= & \max\{\rho((\delta(u))_n,(\delta(v))_n),d(u,v)\}\\
&\le & \max\{\rho(\delta(u),\delta(v)),d(u,v)\}, \quad \text{by (ii)}\\
&=& d(u,v), \quad \text{since } \delta \text{ is nonexpansive}.
\end{eqnarray*}
\item[(iv)] Apply (iii) and the fact that 
$M$ has the AFPP for nonexpansive mappings to get $(z_n)$ in $M$ such that
\[d(z_n, \varphi_n(z_n))\le \frac{1}{n}.\]
The second inequality follows immediately.
\qed\end{enumerate}
\end{proof}

Applying Theorem \ref{quant-BRS} to the family $(T_u)_{u\in M}$ of nonexpansive mappings, we can prove the following result.

\begin{lemma} {\bf (Main technical lemma)}\label{main-lemma}\\
Let $(\lambda_n)$ be a sequence in $[0,1)$, and $K\in\N, \alpha:\N\to\N$ satisfying {\rm (\ref{hyp-K-alpha})}.
Let $b_1,b_2>0$ be such that
\begin{eqnarray}
\forall u\in M\,\exists x^*\in C\big(\rho(\delta(u),x^*)\le b_1
\wedge\rho(x^*,T_u(x^*))\le b_2\big).\label{hyp-main-lemma}
\end{eqnarray}
Then there are sequences $(z_n)$ in $M$ and $(x^*_n)$ in $C$
such that
\[ \forall n\in\N (\rho(\delta(z_n),x^*_n)\le b_1
\wedge \rho \big(x^*_n,T_{z_n}(x^*_n)\big)\le b_2) \]
and for all $\varepsilon >0$,
\[ \forall n\ge 
g(\varepsilon,b_1,b_2,K,\alpha)
\big(d_\infty(((\delta(z_n))_n,z_n), T((\delta(z_n))_n,z_n))\le 
\rho(x^*_n,T_{z_n}(x^*_n))+\varepsilon\big),\]
where
\[ \begin{array}{l} 
g(\varepsilon,b_1,b_2,K,\alpha)=
\max\left\{\left\lceil\frac{1}{\varepsilon}
\right\rceil+1,h(\varepsilon,2b_1+b_2,K,\alpha)\right\},
\end{array}\]
and $h$ is as in Theorem \ref{quant-BRS}.
\end{lemma}
\begin{proof}
First, let us remark that for all $u\in M$, by (\ref{hyp-main-lemma}) and the fact that $T_u$ is nonexpansive, 
there exists an $x^*\in C$ such that 
\begin{eqnarray}
\rho(\delta(u),T_u(\delta(u)))&\le &\rho(\delta(u),x^*)+\rho(x^*, T_u(x^*))+\rho(T_u(x^*),T_u(\delta(u)))\nonumber\\
&\le & 2b_1+b_2.\label{main-lemma-ineq-1}
\end{eqnarray}
Let $(z_n)$ be as in Lemma \ref{useful-lemma}(iv), and  $(x^*_n)$ be the sequence in $C$ obtained by applying  (\ref{hyp-main-lemma}) to $u:=z_n$. 
For every $n\in\N$, $T_{{z_n}}:C\to C$ is nonexpansive, and $(\delta(z_n))_m$ is the Krasnoselski-Mann iteration associated with $T_{z_n}$, starting with $\delta(z_n)\in C$. By the hypothesis and (\ref{main-lemma-ineq-1}), we have that for all $n\geq 0$
\begin{equation}
\rho(\delta(z_n),T_{z_n}(\delta(z_n)))\le 2b_1+b_2 \wedge \rho(\delta(z_n),x^*_n)\le b_1\le 
2b_1+b_2.\label{main-lemma-ineq2}
\end{equation}
Now let $\varepsilon >0$ and $n\ge g(\varepsilon,b_1, b_2,K,\alpha)$. Then $n\ge h(\varepsilon,2b_1+b_2,K,\alpha)$, and we can apply Theorem \ref{quant-BRS} to get that 
\begin{eqnarray*}
\rho\big((\delta(z_n))_n,T_{z_n}((\delta(z_n))_n)\big)\le\rho(x^*_n,T_{z_n}(x^*_n))+\varepsilon.
\end{eqnarray*}
Hence, by Lemma \ref{useful-lemma}(iv)
\begin{eqnarray*}
d_\infty(((\delta(z_n))_n, z_n), T((\delta(z_n))_n, z_n))&\le& \max\left\{\rho(x^*_n, T_{z_n}(x^*_n))+\varepsilon, \frac{1}{n}\right\}\\
&=& \rho(x^*_n, T_{z_n}(x^*_n))+\varepsilon, 
\end{eqnarray*}
since $\displaystyle n\ge g(\varepsilon,b_1,b_2,K,\alpha)\ge \left\lceil\frac{1}{\varepsilon}\right\rceil+1$, so $\displaystyle \frac 1 n<\varepsilon$.
\qed\end{proof}
We are now in position to prove the main theorems of this section.

\begin{theorem}\label{main-rH-rC}
Assume that 
\[\displaystyle\sup_{u\in M}r_C(T_u)<\infty,\]
and $\varphi:\R^*_+\to \R^*_+$ is such that for each $\varepsilon>0$ and $v\in M$ there exists $x^*\in C$ satisfying 
\begin{eqnarray}
\rho(\delta(v),x^*)\le \varphi(\varepsilon)\,\wedge\,\rho(x^*,T_v(x^*))\le 
\sup_{u\in M}r_C(T_u)+\varepsilon.\label{hyp-rH-rC}
\end{eqnarray}
Then
\begin{equation}
r_H(T)\le \sup_{u\in M}r_C(T_u). \label{conclusion-rH-rC}
\end{equation}
\end{theorem}
\begin{proof}
Let $\varepsilon>0$. Define $\displaystyle b_1:=\varphi(\varepsilon) \ \mbox{and} \ b_2:= 
\sup_{u\in M}r_C(T_u)+\varepsilon$. Then, by (\ref{hyp-rH-rC}), we get that 
\begin{equation*}
\forall v\in M\exists x^*\in C\big(\rho(\delta(v),x^*)\le b_1
\wedge\rho(x^*,T_v(x^*))\le b_2\big),
\end{equation*}
so (\ref{hyp-main-lemma}) is satisfied.  
Let $(\lambda_n), K,\alpha$ be as in the hypothesis of 
Lemma \ref{main-lemma}.  Then we can apply that lemma to get sequences 
$(z_n)$ in $M$ and $(x^*_n)$ in $ C$  such that 
\[ \forall n\in\N (\rho(x^*_n,T_{z_n}(x^*_n))\le b_2), \] 
and for all $n\ge g(\varepsilon,b_1,b_2,K,\alpha)$,
\[d_\infty(((\delta(z_n))_n,z_n), T((\delta(z_n))_n,z_n))\le 
\rho(x^*_n,T_{z_n}(x^*_n))+\varepsilon,\]
where $g(\varepsilon,b_1,b_2,K,\alpha)$ is defined as in Lemma \ref{main-lemma}. 
It follows that
\begin{eqnarray*} 
 r_H(T)&\le & d_\infty(((\delta(z_n))_n,z_n), T((\delta(z_n))_n,z_n))\le\rho(x^*_n,T_{z_n}(x^*_n))+\varepsilon\\
 &\le & b_2+\varepsilon=\sup_{u\in M}r_C(T_u)+2\varepsilon.
\end {eqnarray*}
Since $\varepsilon>0$ was arbitrary, (\ref{conclusion-rH-rC}) follows.
\qed\end{proof}

The next theorem is obtained by applying Theorem \ref{quant-Ishikawa} to the family $(T_u)_{u\in M}$.

\begin{theorem}\label{main-x-yn-bounded}
Let $(\lambda_n), K,\alpha$ be as in the hypothesis of Theorem \ref{quant-BRS}.\\
Assume that there is $b>0$ such that 
\begin{eqnarray}
\forall u\in M\exists y\in C\big(\rho(\delta(u),y)\le b\wedge \forall m,p\in\N (\rho(y^u_m,y^u_p)\le b)\big), \label{hyp-main-x-yn-bounded}
\end{eqnarray}
where $(y^u_n)$ is the Krasnoselski-Mann iteration from $y$, associated with $T_u$.

Then $r_H(T)=0$. 
\end{theorem}
\begin{proof}
Take $(z_n)$ as in Lemma \ref{useful-lemma}(iv).  By hypothesis, for every $n\in\N$ there exists $y\in C$ such that
\[\displaystyle \rho(\delta(z_n),y)\le b\wedge \forall m,p\in\N\, (\rho(y_m^{z_n},y_p^{z_n})\le b),\]
Let now $\varepsilon >0$ and $n\ge \tilde{g}(\varepsilon,b,K,\alpha)$.
We are in position to apply Theorem \ref{quant-Ishikawa} to get that 
\[\rho\big((\delta(z_n))_n,T_{z_n}((\delta(z_n))_n)\big)\le\varepsilon.\]
Then, by Lemma \ref{useful-lemma}(iv)
\begin{eqnarray*}
d_\infty(((\delta(z_n))_n, z_n), T((\delta(z_n))_n, z_n))&\le& \max\left\{\varepsilon, \frac{1}{n}\right\}\\
&=& \varepsilon, 
\end{eqnarray*}
since $\displaystyle n\ge \tilde{g}(\varepsilon,b,K,\alpha)\ge \left\lceil\frac{1}{\varepsilon}\right\rceil+1$, so $\displaystyle \frac 1 n<\varepsilon$.
In particular, it follows that $r_H(T)=0$.
\qed\end{proof}

The above theorems have some straightforward consequences.

\begin{corollary}\label{cor-main-rH-rC-used}
Assume that $\varphi:\R^*_+\to \R^*_+$ is such that
\begin{equation}
\forall\varepsilon>0\forall u\in M\exists x^*\in C\big(\rho(\delta(u),x^*)\le 
\varphi(\varepsilon)\wedge\rho(x^*,T_u(x^*))\le\varepsilon\big).\label{hyp-C-AFPP}
\end{equation}
Then $r_H(T)=0$.
\end{corollary}
\begin{proof}
From the hypothesis, it follows immediately that $r_C(T_u)=0$ for all $u\in M$. Thus, $\displaystyle\sup_{u\in M}r_c(T_u)=0$. Apply now Theorem \ref{main-rH-rC}.
\qed\end{proof}

\begin{corollary}\label{cor-useful}
Assume that for all $u\in M$, the Krasnoselski-Mann iteration $(\delta(u))_n$ is bounded. Then $r_H(T)=0$.   
\end{corollary}
\begin{proof}
Apply Theorem \ref{main-x-yn-bounded} with $y:=\delta(u)$.
\qed\end{proof}

Theorem \ref{th-kirk-CAT} is an immediate consequence of Corollary \ref{cor-useful}, since if $C$ is bounded, then $(\delta(u))_n$ is bounded by the diameter of $C$ 
for each $u\in M$. 

\subsection{Families of unbounded convex sets}

\noindent In this section we indicate that 
all the above results can be generalized to families $(C_u)_{u\in M}$ of nonempty {\em unbounded} convex subsets of the hyperbolic space $(X,\rho,W)$.

Let $(C_u)_{u\in M}\subseteq X$ be a family of convex sets such that there exists a nonexpansive {\em selection } function $\delta:M\to \bigcup_{u\in M}C_u$, that is a nonexpansive function with the property  
\begin{equation} \forall u\in M\big(\delta(u)\in C_u\big). \label{selection}
\end{equation}

We consider the following subspace of $(X\times M)_\infty$:
\[ H:=\{(x,u): u\in M, x\in C_u\}.\]
Let $P_1:H\to \bigcup_{u\in M}C_u, P_2:H\to M$ be the projections. In the sequel, we consider nonexpansive functions  $T:H\to H$ satisfying
\begin{equation}
 \forall (x,u)\in H\,\, \big((P_1\circ T)(x,u)\in C_u\big).\label{hyp-T-Cu}
\end{equation}
It is easy to see that we can define for each $u\in M$ the nonexpansive function
\[T_u:C_u\to C_u, \quad T_u(x)=(P_1\circ T)(x,u).\]
For each $u\in M$, we denote the Krasnoselski-Mann iteration starting from $x\in C_u$ and associated with $T_u$ by $(x^u_n)$ (or by $(x_n)$, when $u$ is clear from the context).

For each $n\ge 0$, we define 
\begin{eqnarray*}
\varphi_n:M\to M, & \varphi_n(u)=(P_2\circ T)((\delta(u))_n,u).
\end{eqnarray*}
The following results can be proved in a similar manner as Lemmas \ref{useful-lemma}, \ref{main-lemma}, and  Theorems \ref{main-rH-rC}, \ref{main-x-yn-bounded}.

\begin{lemma}\label{useful-lemma-family}
\begin{enumerate}
\item[(i)] for all $n\ge 0,\, u,v\in M$, and $x\in C_u,y\in C_v$
\begin{eqnarray}
\rho(x_n^u, y_n^v)\le \max\{\rho(x,y),d(u,v)\}=d_\infty((x,u),(y,v));
\end {eqnarray}
\item[(ii)] for all $n\ge 0$, $\varphi_n$ is nonexpansive;
\item[(iii)] there is a sequence $(z_n)$ in $M$ such that for all $n\ge 0$,
\[\begin{array}{l}
\displaystyle d(z_n, \varphi_n(z_n))\le \frac{1}{n}, \text{and}\\
\displaystyle d_\infty(((\delta(z_n))_n, {z_n}), T((\delta(z_n))_n, z_n))\le \max\left\{\rho((\delta(z_n))_n, T_{z_n}(((\delta(z_n))_n, \frac{1}{n}\right\}.
\end{array}
\]
\end{enumerate}
\end{lemma}

\begin{lemma}\label{main-lemma-family}
Let $(\lambda_n), K,\alpha$ be as in the hypothesis of Theorem \ref{quant-BRS}.\\
Assume that $b_1,b_2>0$ are such that
\begin{eqnarray}
\forall u\in M\,\exists x^*\in C_u\big(\rho(\delta(u),x^*)\le b_1
\wedge\rho(x^*,T_u(x^*))\le b_2\big).
\end{eqnarray}
Then there are sequences $(z_n)$ and $(x^*_n)$ such that
\[ \forall n\in\N (z_n\in M\wedge x^*_n\in C_{z_n}\wedge\rho(\delta(z_n),x^*_n)\le b_1
\wedge \rho \big(x^*_n,T_{z_n}(x^*_n)\big)\le b_2) \]
and for all $\varepsilon >0$,
\[ \forall n\ge 
g(\varepsilon,b_1,b_2,K,\alpha)
\big(d_\infty(((\delta(z_n))_n,z_n), T((\delta(z_n))_n,z_n))\le 
\rho(x^*_n,T_{z_n}(x^*_n))+\varepsilon\big),\]
where $g(\varepsilon,b_1,b_2,K,\alpha)$ is as in Lemma \ref{main-lemma}.
\end{lemma}

\begin{theorem}\label{main-family-rH-rC}
Assume that 
\[\sup_{u\in M}r_{C_u}(T_u)<\infty.\]
Suppose also that $\varphi:\R^*_+\to \R^*_+$ is such that for each $\varepsilon>0$ and $v\in M$ there exists $x^*\in C_v$ satisfying 
\begin{eqnarray}
\rho(\delta(v),x^*)\le \varphi(\varepsilon)\,\wedge\,\rho(x^*,T_v(x^*))\le 
\sup_{u\in M}r_C(T_u)+\varepsilon.\label{hyp-main-family-rH-rC}
\end{eqnarray}
Then
\[r_H(T)\le \sup_{u\in M}r_{C_u}(T_u).\]
\end{theorem}

\begin{theorem}\label{main-family-x-yn-bounded}
Let $(\lambda_n), K,\alpha$ be as in the hypothesis of Theorem \ref{quant-BRS}.\\
Assume that there is $b>0$ such that 
\begin{eqnarray}
\forall u\in M\exists y\in C_u\big(\rho(\delta(u),y)\le b\wedge \forall m,p\in\N (\rho(y^u_m,y^u_p)\le b)\big).\label{Cu-hyp-th-3.21}
\end{eqnarray}
Then $r_H(T)=0$.
\end{theorem}

We get also the following corollary:

\begin{corollary}
Assume that $(C_u)_{u\in M}\subseteq X$ is a family of nonempty {\bf bounded} convex sets such that there is $b>0$ with the property that 
\begin{equation}
\forall u\in M \big(diam(C_u)\le b\big). \label{hyp-family-bounded}
\end{equation}
Then $H$ has the AFPP for nonexpansive mappings $T\!:\!H\!\to\!H$ satisfying (\ref{hyp-T-Cu}).
\end{corollary}
\begin{proof} 
The hypothesis of Theorem \ref{main-family-x-yn-bounded} is satisfied with $y:=\delta(u)$.
%Let $T:H\to H$ be a nonexpansive function satisfying (\ref{hyp-T-Cu}). Since for all $u\in M$,  $C_u$ is bounded, we get that $r_{C_u}(T_u)=0$, so $\sup_{u\in M}r_{C_u}(T_u)=0$. Then, if we define $\varphi(\varepsilon)=b$, 
%we get (\ref{hyp-main-family-rH-rC}), so we can  apply Theorem \ref{main-family-rH-rC}.
\qed\end{proof}

\section{Uniform approximate fixed point property}

Let $(X,\rho)$ be a metric space, $C\subseteq X$ a nonempty subset, and $\mathcal{F}$ be a class of functions $T:C\to C$. We say that $C$ has the {\em uniform approximate fixed point property (UAFPP)} for $\mathcal{F}$ if for all $\varepsilon>0$ and $b>0$ there exists an $D>0$ such that for 
each point $x\in C$ and for each function $T\in{\mathcal F}$,\begin{equation}
\rho(x,T(x))\le b\Rightarrow \exists x^*\in C\big(\rho(x,x^*)\le D\wedge 
\rho(x^*, T(x^*))\leq\varepsilon\big).
\label{uafpp-def}
\end{equation}

Let $(X,\rho,W)$ be a hyperbolic space, and $C\subseteq X$ be a nonempty convex subset. Assume that $(\lambda_n)$ is a sequence in $[0,1]$. We say that $C$ has the 
{\em $\lambda_n$-uniform asymptotic regularity property } if for all $\varepsilon>0$ 
and $b>0$ there exists an $N\in\N$ such that for each point $x\in C$ and for 
each nonexpansive function $T:C\to C$, 
\begin{equation}
\rho(x,T(x))\le b\Rightarrow \forall n\ge N\big(\rho(x_n,T(x_n))\leq\varepsilon\big), \label{lambda-uar-def}
\end{equation}
where $(x_n)$ is the Krasnoselski-Mann iteration defined as in (\ref{KM-lambda-n-def}).

These definitions are inspired by Theorem \ref{quant-BRS}, our quantitative version of Bor-wein-Reich-Shafrir Theorem, which can be used to prove the following equivalent characterizations. 

\begin{proposition}
Let $(X,\rho,W)$ be a hyperbolic space, and $C\subseteq X$ be a nonempty convex subset. Then the following are equivalent:
\begin{enumerate}
\item[(i)] $C$ has the UAFPP for nonexpansive functions;
\item[(ii)] $C$ has the $\lambda_n$-uniform asymptotic regularity property for 
some $(\lambda_n)$  in $[0,1]$; 
\item[(iii)] $C$ has the $\lambda_n$-uniform asymptotic regularity property for all 
$(\lambda_n)$  in $[0,1]$ which are divergent in sum and bounded away from 1.
\end{enumerate}
\end{proposition}
\begin{proof}
$(iii)\Rightarrow (ii)$ is obvious.\\
$(i)\Rightarrow (iii)$ Let $\varepsilon>0, b>0$, and $D>0$ be such that (\ref{uafpp-def}) holds with ${\mathcal F}$ being the class of nonexpansive functions. Let $(\lambda_n)$  in $[0,1]$ be divergent in sum and bounded away from 1. Then  there are $K\in\N$, and $\alpha:\N\to\N$ such that (\ref{hyp-K-alpha}) holds.  Let $x\in C$ and $T:C\to C$ nonexpansive 
be such that $\rho(x,T(x))\le b$. Then there is $x^*\in C$ such that 
$\rho(x,x^*)\le D,$ and $\rho(x^*, T(x^*))\leq\varepsilon$. Let $b':=\max\{b,D\}$. Then
\[\rho(x,T(x))\le b', \,\,\rho(x,x^*)\le b',\]
so we can apply Theorem \ref{quant-BRS} to get 
$N:=h(\varepsilon,b',K,\alpha)$ such that
\[\forall n\ge N\big(\rho(x_n,T(x_n))\leq \rho(x^*,T(x^*))+\varepsilon\leq 2\varepsilon\big).\]
$(ii)\Rightarrow(i)$ Let $(\lambda_n)$ in $[0,1]$ be such that $C$ has the
$\lambda_n$-uniform asymptotic regularity property. Let $\varepsilon>0, b>0$, 
and $N\in\N$ as in  (\ref{lambda-uar-def}). Let now $x\in C$ and $T:C\to C$ 
nonexpansive be such that $\rho(x,T(x))\le b$. Then $\rho(x_N, T(x_N))\leq\varepsilon$, and 
\begin{eqnarray*}
\rho(x,x_N) & \le & \sum_{i=0}^{N-1}\rho(x_i,x_{i+1})=
\sum_{i=0}^{N-1}\lambda_i\rho(x_i,T(x_i))\\
& \le & \rho(x,T(x))\sum_{i=0}^{N-1}\lambda_i, \text{ since } 
(\rho(x_n,T(x_n))) \text{ is nonincreasing}\\
&\le & b\sum_{i=0}^{N-1}\lambda_i.
\end {eqnarray*}
Let us take now $\displaystyle D:= b\sum_{i=0}^{N-1}\lambda_i$, and $x^*:=x_N$. 
Then (\ref{uafpp-def}) holds.
\qed\end{proof}

Let us remark the following fact. A first attempt to define  the property 
that $C$ has the uniform approximate fixed point property for nonexpansive 
functions is in the line of the Goebel-Kirk Theorem \cite{Goebel/Kirk}, that is  
\begin{eqnarray}
 \forall \varepsilon>0\,\exists D>0\,\forall x\in C\,\forall\, T:C\to C 
\text{ nonexpansive }\hfill\nonumber \\
\exists x^*\in C\big(\rho(x,x^*)\le D\wedge \rho(x^*, T(x^*))\leq\varepsilon\big). 
\label{uafpp-goebel-kirk}
\end {eqnarray}

If $C$ is bounded, then $C$ satisfies (\ref{uafpp-goebel-kirk}) 
by \cite{Goebel/Kirk}.  But also conversely, if $C$ 
satisfies (\ref{uafpp-goebel-kirk}) (even if only for 
constant functions $T$) then $C$ is already bounded: \\ 
Assume that $C$ satisfies (\ref{uafpp-goebel-kirk}) for all 
constant functions $T$. 
Then for $\varepsilon:=1$ we get $D_1\in\N$ such that for all $x\in C$, and for 
all constant functions $T:C\to C$, there is $x^*\in C$ with $\rho(x,x^*)\le D_1$ and 
$\rho(x^*, T(x^*))\leq 1$. It follows that 
\begin{equation}
\rho(x,T(x))\le \rho(x,x^*)+\rho(x^*, T(x^*))+\rho(T(x^*),T(x))
\le 2D_1+1\label{ineq-unif-GK}
\end{equation}
Now, if we assume that $C$ is unbounded, there are $x,y\in C$ such that 
$\rho(x,y)>2D_1+1$. Define $T:C\to C, \, T(z)=y$ for all $z\in C$. Then 
$\rho(x,T(x))=\rho(x,y)>2D_1+1$ which contradicts (\ref{ineq-unif-GK}).

The next results give some partial answers to problem 27 in \cite{Kirk(04)} 
which asks whether the product $H=(C\times M)_\infty$ of a closed convex subset $C$ of a 
complete CAT(0)-space $X$ (having the 
geodesic line extension property) with the AFPP (or -- equivalently -- being 
geodesically bounded) and a metric space $M$ having the AFPP again has 
the AFPP: we show that this is true if $C$ has the UAFPP (even in the case where 
$X$ is just a hyperbolic space) and a technical condition is satisfied 
which, in particular, holds if $M$ is bounded:

\begin{theorem}\label{theorem-UAFPP+AFPP}
Let $(X,\rho,W)$ be a hyperbolic space, $C\subseteq X$ a convex subset with the UAFPP for nonexpansive functions, $(M,d)$ be a metric space with the AFPP for nonexpansive functions, and $\delta:M\to C$ be a nonexpansive selection function.
Let $T:H\to H$ be a nonexpansive function such that there exists $b>0$ satisfying
\begin{equation} 
\forall u\in M \big(\rho(T_u(\delta(u)),\delta(u))\le b\big). \label{hyp-theorem-Kirk-open}
\end{equation}
Then $r_H(T)=0$.
\end{theorem}
\begin{proof} 
Let $\varepsilon>0$. Since $C$ has the UAFPP for nonexpansive functions, there exists $D>0$ (depending on $\varepsilon$ and $b$) such that (\ref{uafpp-def}) holds for each nonexpansive function. For each $u\in M$, using (\ref{hyp-theorem-Kirk-open}) and  applying (\ref{uafpp-def}) for $x:=\delta(u)$ and $T_u$, we get $x^*\in C$ such that $\rho(\delta(u),x^*)\le D$ and $\rho(x^*, T_u(x^*))\leq\varepsilon $. 
Hence, the hypothesis of Corollary \ref{cor-main-rH-rC-used} is satisfied with $\varphi(\varepsilon):=D$, so $r_H(T)=0$ follows.
\qed\end{proof}

\begin{corollary}
Let $(X,\rho,W)$ be a hyperbolic space, and $C\subseteq X$ a convex subset with the UAFPP for nonexpansive functions.
Assume that $(M,d)$ is a {\bf bounded} metric space with the AFPP for nonexpansive functions.

Then  $H:=(C\times M)_\infty$ has the AFPP for nonexpansive functions. 
\end{corollary}
\begin{proof}
Let $x\in C$ be arbitrary, and define $\delta:M\to C$ by $\delta(u)=x$. Let $T:H\to H$ be a nonexpansive function. Fix some $u_0\in M$, and define $b:=\rho(x,T_{u_0}(x))+diam(M)$. Then for each $u\in M$, we have that 
\begin{eqnarray*}
\rho(x,T_u(x))&\le& \rho(x,T_{u_0}(x))+ \rho(T_{u_0}(x), T_u(x))\\
&\le & \rho(x,T_{u_0}(x))+ d_\infty(T(x,u_0), T(x,u))\\
&\le &\rho(x,T_{u_0}(x))+ d_\infty((x,u_0), (x,u))\\
&= &\rho(x,T_{u_0}(x))+ d(u_0,u)\le b.
\end {eqnarray*}
Thus, we can apply Theorem \ref{theorem-UAFPP+AFPP} to conclude that  $r_H(T)=0$. 
\qed\end{proof}

We end this section by remarking that, using the same ideas, 
we can define the notion of $C$ having the uniform fixed point property. 
Thus, $C$ has the {\em uniform fixed point property (UFPP)} for $\mathcal{F}$ 
if for all $b>0$ there exists $D>0$ such that for each point $x\in C$ and 
for each function $T\in{\mathcal F}$,\begin{equation}
\rho(x,T(x))\le b\Rightarrow \exists x^*\in C\big(\rho(x,x^*)\le D\wedge T(x^*)=x^*\big).
\label{ufpp-def}
\end{equation}

\begin{proposition}
Assume that $(X,\rho)$ is a complete metric space. 
Let $\mathcal{F}$ be the class of contractions with a common contraction constant 
$k\in (0,1)$. Then each closed subset $C$ of $X$ has the UFPP for $\mathcal{F}$.
\end{proposition}
\begin{proof}
By Banach's Contraction Mapping Principle we know that each mapping $T\in\mathcal{F}$ has a unique fixed point $x_0$, and, moreover, for each $x\in C$, 
\begin{equation}
\rho(T^n(x), x_0)\le \frac{k^n}{1-k}\rho(x,T(x)) \quad \text{ for all } n\ge 0.\label{contraction-Banach-estimation}
\end{equation}
For $n=0$, this yields $\displaystyle\rho(x,x_0)\le \frac{1}{1-k}\rho(x,T(x))$. Then 
$\displaystyle \rho(x,T(x))\le b$ implies $\displaystyle\rho(x,x_0)\le \frac{1}{1-k}b.$ Hence, (\ref{uafpp-def}) holds with $\displaystyle D:=\frac{1}{1-k}b$ and $x^*:=x_0$.
\qed\end{proof}
It is known that all directionally bounded closed convex subsets of complete hyperbolic 
spaces in the sense of \cite{Reich/Shafrir} have the approximate fixed 
point property for all nonexpansive mappings and that there are 
(even for normed spaces) unbounded but directionally bounded convex 
subsets (\cite{Shafrir(91)}). We conclude this paper with an 

{\bf Open Problem:} Are there  unbounded convex subsets of some 
hyperbolic space which have the UAFPP for all nonexpansive mappings 
$T:C\to C?$

\end{document}